\documentclass[12pt]{amsart}
\usepackage[all]{xy}
\usepackage{amsmath,amssymb}
\theoremstyle{plain}

\newtheorem{theorem}{Theorem} [section]
\newtheorem{proposition}[theorem]{Proposition}
\newtheorem{theoremaa}{Theorem}
\newtheorem{proposition.a}[theoremaa]{Proposition}
\newtheorem{corollary.a}[theoremaa]{Corollary}
\newtheorem{corollary}[theorem]{Corollary}
\newtheorem{lemma}[theorem]{Lemma}
\theoremstyle{definition}
\newtheorem{definition}[theorem]{Definition}
\newtheorem{remark}[theorem]{Remark}
\newtheorem{example}[theorem]{Example}

\numberwithin{equation}{section}

\theoremstyle{remark}

\newcommand{\lra}{\longrightarrow}
\newcommand{\noi}{\noindent}
\newcommand{\PP}{\mathbf{P}}

\newcommand{\RR}{\mathbf{R}}
\newcommand{\ZZ}{\mathbf{Z}}
\newcommand{\NN}{\mathbf{N}}
\newcommand{\CC}{\mathbf{C}}
\newcommand{\QQ}{\mathbf{Q}}
\newcommand{\AAA}{\mathbf{A}}
\newcommand{\OO}{\mathcal{O}}
\newcommand{\II}{\mathcal{I}}
\newcommand{\QQQ}{\mathcal{Q}}

\newcommand{\FF}{\mathcal{F}}
\newcommand{\fra}{\mathfrak a}

\newcommand{\frq}{\mathfrak q}

\newcommand{\JJ}{\mathcal{J}}

\newcommand{\HH}[3]{H^{{#1}} \big( {#2} , {#3} \big) }
\newcommand{\hh}[3]{h^{{#1}} \big( {#2} , {#3} \big) }
\newcommand{\HHH}[3]{H^{{#1}} \Big(\, {#2} \, , \,  {#3} \, \Big) }

\newcommand{\for}{ \ \ \text{ for } \ }
\newcommand{\fall}{ \ \ \text{ for all } \ }

\newcommand{\reg}{\text{reg}}

\newcommand{\pr}{\prime}

\newcommand{\Bl}{\text{Bl}}
\newcommand{\adeg}{\text{adeg}}

\begin{document}

\title {Positivity and Complexity of Ideal Sheaves}

\author{Steven Dale Cutkosky}
\address{Department of Mathematics \\ University of
Missouri \\ Columbia, MO  }
\thanks{Research of first author partially supported
by  NSF Grant DMS 99-88533.}
\email{cutkoskys@missouri.edu}
\author{Lawrence Ein} 
\address{Department of Mathematics \\
University of Illinois at Chicago \\ 851 South Morgan
St., M/C. 249
\\ Chicago, IL  60607-7045}
\thanks{Research of second author  partially
supported by NSF Grant DMS 99-70295}
\email{ein@math.uic.edu}
\author{Robert Lazarsfeld} 
\address{Department of Mathematics
\\ University of Michigan \\ Ann Arbor, MI  48109}
\email{rlaz@math.lsa.umich.edu}
\thanks{Research of  third author  partially
supported by the J. S. Guggenheim Foundation and NSF
Grant DMS 97-13149}

\maketitle
\setcounter{section}{0}

\section*{Introduction} 
The problem of bounding the
``complexity" of a polynomial ideal in terms of the
degrees of its generators has attracted a great deal of
interest in recent years. Results in this direction go
back at least as far as the classical work \cite{Hermann} of 
Hermann on the ideal membership problem, and the
effective Nullstellensatz of Brownawell
\cite{Brownawell} and Koll\'ar \cite{Kollar} marks a major recent
advance. With the development of computational
algebraic geometry the question has taken  on
increasing importance, and  it came into particularly
clear focus through the influential paper  \cite{BM}
of Bayer and Mumford. More recently, the theorem of
\cite{CHT} and \cite{Kod} concerning regularity
of powers raises the question of bounding the
complexity of powers of an ideal, and suggests that
asymptotically the picture   should
become very clean. 

The aim of the present paper is to examine some of the
results and questions of
\cite{BM}, \cite{Vas} and \cite{Cut.Irr} 
from a geometric perspective, in the spirit of
\cite{EL99}. Our thesis is that much of this  material
is clarified, and parts rendered
 transparent, when viewed through the
lenses of vanishing theorems and intersection theory.
 Specifically, motivated by the work
\cite{Paoletti.JDG}, \cite{Paoletti2} of Paoletti we introduce an
invariant $s(\JJ)$ that measures in effect how much
one has to twist an ideal
$\JJ$ in order to make it  positive.
Degree bounds on   generators of $\JJ$ yield
bounds on this  $s$-invariant,   but in
general
$s(\JJ)$ can be small even when the degrees of
generators are large. We  prove that the 
$s$-invariant $s(\JJ)$ computes the asymptotic
regularity of large powers of an ideal sheaf, and
bounds the asymptotic behaviour of several other
natural measures of complexity considered in
\cite{BM} and \cite{Vas}. We also show that this invariant behaves
very well with respect to natural geometric and
algebraic operations.
 This  leads for example to a
considerably  simplified analogue of the
construction from 
\cite{Cut.Irr} of varieties with irrational asymptotic
regularity.

Turning to  a more detailed overview of the contents
of this paper, we start by fixing the setting in which
we shall work. Denote by $X$ an irreducible
non-singular projective variety of dimension $n$
defined over an uncountable algebraically closed field
$K$ of arbitrary characteristic, and let
$H$ be a fixed ample divisor (class) on
$X$.  The most natural and important example is of
course $X = \PP^n$ and $H$ a hyperplane, and in fact
little essential will be lost to the reader who
focuses on this classical case. However since we are
working geometrically it seems natural  to consider
general varieties, and in the end it is no harder to do
so. Given an ideal sheaf $\JJ
\subseteq \OO_X$,  consider the blowing-up
$\nu : 
\Bl_{\JJ}(X)
\lra X$   of $\JJ$, with exceptional
divisor
$F$. We define
\[
s_H(\JJ) \  = \ \inf \{s \in \RR \mid  \nu^* (sH) - F
\ \text{is an ample $\RR$-divisor on $\Bl_{\JJ}(X) $}
\, \}.\]
One has $s_H(\JJ) \le d_H(\JJ)$, where 
$d_H(\JJ)$ denotes the least integer $d > 0 $ such
that $\JJ(dH)$ is globally generated, but in general
the inequality is strict.\footnote{In fact, it is
 possible for $s_H(\JJ)$ to be irrational.} This
\textit{s-invariant} is closely related the Seshadri
constants introduced by Demailly, and has been studied
by Paoletti when $\JJ$ is the ideal sheaf of a smooth
subvariety of $X$. (See Remark \ref{S.Const.Rem} below). In a
general way, our goal is to bound the  ``complexity" 
of $\JJ$ (or at least its powers) in terms of this
invariant. 

We do not use the term ``complexity" here in any
technical sense. Rather, guided by
\cite{BM}, we  consider various natural invariants
that each give a picture of how complicated one might
consider
$\JJ$ to be:
\begin{itemize}
\item[(0).] The degrees of the irreducible components
of the zero-locus
$\text{Zeroes}(\JJ) \subseteq X$ of $\JJ$;
\vskip 3pt
\item[(1).] The degrees of all the ``associated
subvarieties" of
$\text{Zeroes}(\JJ)$, including those corresponding to
the embedded primes in a primary
decomposition of
$\JJ$; 
\vskip 3pt
\item[(2).] For $H$ very ample, the
\textit{Castelnuovo-Mumford regularity} 
$\text{reg}(\JJ)$ of $\JJ$, which   measures roughly
speaking the cohomological complexity of $\JJ$;
\vskip 3pt
\item[(3).] The \textit{index of nilpotency}
$\text{nilp}(\JJ)$ of
$\JJ$, i.e. the least integer $t > 0$ such that 
\[ \Big( \sqrt{\JJ} \Big) ^t \subseteq \JJ. \]
\end{itemize}
 In settings (0) and (1), one can ask
also  for degree bounds after having attached
multiplicities to the components in question: allowing
embedded components, this leads to what Bayer and
Mumford call the \textit{arithmetic degree} of $\JJ$.
The index of nilpotency is closely related to the
effective Nullstellensatze of  Koll\'ar and \cite{EL99},
and various relations among these invariants have been
established (\cite{BM}, \cite{Vas},
\cite{MV}). In the classical situation, where
$\JJ$ is replaced by a  homogeneous ideal 
$I$ generated by forms of degree $d$, it is elementary
to obtain a Bezout-type bound for (0), while the main
theorem of  \cite{Kollar} gives the analogous statement for
(3).  However,   Bayer and Mumford observe that there
cannot exist singly exponential bounds in $d$ for the
regularity or arithmetic degree.

In the direction of (0) and (1), one has: 
\begin{proposition.a} Let $s = s_H(\JJ)$. Then
\begin{equation} \sum \, s^{\dim Z} \cdot \deg_H Z \
\le 
\ s^n \cdot  \deg_H X,  
\end{equation} where the sum is taken over all
irreducible components of $\text{Zeroes}(\JJ)$. If
$\JJ$ is integrally closed, then the same inequality
holds including in the sum also the embedded
associated subvarieties of
$\text{Zeroes}(\JJ)$.
\end{proposition.a}
 \noi The first assertion  follows
already from the positivity theorems for
 intersection classes established in \cite{FL},
although the elementary direct approach of
\cite{EL99} also applies. The stronger statement
for integrally closed ideals, while elementary, seems
to have been overlooked. We also give examples
(Example \ref{Pathology.Example}) to show that if
$\JJ$ is not integrally closed, then one cannot bound
the number of embedded components in terms of
$s_H(\JJ)$.  

Assume now that $H$ is very ample. In this setting
the Castelnuovo-Mumford regularity  $\reg_H(\JJ)$ of
$\JJ$ (with respect to $\OO_X(H)$) can be defined just
as in the classical case $X = \PP^n$. 
\begin{theoremaa} \label{Asympt.Reg.Thm} One has the
equalities 
\[ \lim_{p \to \infty} \frac{ \reg_H(\JJ^p)}{p}   \ =
\lim_{p \to \infty} \frac{ d_H(\JJ^p)}{p} \ = \ 
s_H(\JJ),
\]
where given an ideal sheaf $\II \subseteq \OO_X$,
$d_H(\II)$ denotes as above the least integer $d$ such
that
$\II(dH)$ is globally generated. 
\end{theoremaa}
\noi  Thus the $s$-invariant $s_H(\JJ)$ governs
exactly the asymoptotic regularity of powers of
$\JJ$. As indicated above, this  result was suggested
by  the theorems  of 
\cite{CHT} and \cite{Kod}, which
prove  the analogue of the first equality for
homogeneous ideals. 

Continuing to assume that $H$ is very ample, 
one can define as in
\cite{BM} the codimension $k$ contribution
$\text{adeg}_H^k(\JJ)$ to arithmetic degree of
$\JJ$, which measures (taking into account suitable
multiplicities) the degrees of the
codimension $k$ irreducible and embedded components of
the scheme defined by
$\JJ$. As in
\cite{BM} there are upper bounds --- at least
asymtotically --- for this degree in terms of the
regularity, and we deduce
\begin{corollary.a} Denote by $\text{adeg}_H^k(\JJ^p)$
the codimension $k$ contribution to the arithmetic
degree of $\JJ^p$. Then 
\[  
\limsup_{p \to \infty}
\frac{\text{adeg}_H^k(\JJ^p)}{p^k}
\ \le \ \frac{s_H(\JJ)^k}{k!} \cdot \deg_H(X). 
\]
\end{corollary.a}
\noi In general this statement is the best
possible: for instance equality holds for complete
intersections of hypersurfaces of the same degree in
projective space. As in the case of Theorem B, the
simple asymptotic statement contrasts with the
examples presented in
\cite{BM} showing that there cannot be a singly
exponential bound for $\text{adeg}^k(\JJ)$ in terms of
$d_H(\JJ)$ (let alone in terms of $s_H(\JJ)$). 

Turning finally to the index of nilpotency,
one can canonically attach to
$\JJ$ an integer $r(\JJ)$ arising as the maximum of
the multiplicities of the irreducible
components of the exceptional divisor of the
normalized blow-up of $\JJ$.
These multiplicities appear  in a
Bezout-type bound strengthening Proposition A, which
in particular gives rise
 to
the inequality
$r(\JJ) \le  s_+^n 
\cdot \deg_H X  $, where $s_+ = \max \{ 1, s_H(\JJ)
\}$.  The results of
\cite{EL99} (and also \cite{Hickle}) show that $ \big
( \sqrt{\JJ} \big ) ^{ n \cdot r(\JJ)  }
 \subseteq \JJ$ and more generally that
\begin{equation} \left ( \sqrt{\JJ} \right ) ^{  r(\JJ)
\cdot (n + p - 1)   }
\ \subseteq   \JJ^{\, p} \tag{*} \end{equation}
for all $p \ge 1$. Note that (*) leads to the
  asymptotic statement 
$ \limsup_{p \to \infty} \frac{\left(
\, \text{nilp}(\JJ^p)  \,
\right)  }{p}    \ \le  s_+^n \cdot \deg_H X $. 

Motivated by the influence of the $s$-invariant
 in these questions, we study its behavior
under natural geometric and algebraic operations. In
this direction we prove for example:
\begin{proposition.a} \label{alg.properties.s.inv}
Let $\JJ_1, \JJ_2 \subseteq \OO_X$ be ideal sheaves on
$X$. Then one has:
 \begin{gather*}
s_H\big(\JJ_1 \cdot \JJ_2\big) \le s_H(\JJ_1) +
s_H(\JJ_2) \\
s_H \big( \JJ_1 + \JJ_2 \big) \le \max \big \{
s_H(\JJ_1) \, , \, s_H(\JJ_2) \, \big \}.
\end{gather*}
Moreover, if $\overline \JJ$ denotes the integral
closure of an ideal $\JJ$, then $s_H(\overline \JJ) =
s_H(\JJ)$.
\end{proposition.a}
\noi In view of Theorem \ref{Asympt.Reg.Thm}, this
result shows that the asymptotic regularity of large
powers of an ideal satisfies much better formal
algebraic properties than are known or expected to
hold for the regularity of an ideal itself. 

Our exposition is organized as follows. In \S1 we
define and study the $s$-invariant measuring the
positivity of an ideal sheaf.  Degree and nilpotency
bounds --- which for the most part involve only minor
modifications to results from
\cite{EL99} --- are given in \S2.
Finally, in \S 3 we consider asymptotic bounds on the
regularity and arithmetic degree of large powers of
an ideal. 

\vskip 10pt

\section{The $s$-invariant of an ideal sheaf}

In the present section we define and study the
$s$-invariant of an ideal sheaf with respect to an
ample divisor. 

 We start by fixing some notation that will remain in
force throughout the paper. Let
$X$ be a non-singular irreducible projective variety
defined over an uncounctble algebraically closed field
$K$ of arbitrary characteristic, and consider  a
fixed  coherent ideal sheaf
$\JJ \subseteq
\OO_X$.   Denote by
\[ \nu  : W =
\text{ Bl}_{\JJ}(X) \lra X \] the blowing up of $X$
along $\JJ$.  Then $\JJ$ becomes locally principal on
$W$, i.e there is an 
effective Cartier divisor  $F$ on $W$ (namely the
excptional divisor of $\nu$) such that
\[ \JJ \cdot \OO_{W} = \OO_{W}(-F) .\] 

Fix now an ample divisor (class) $H$ on $X$. If $s \gg
0$, then
$\nu^*(sH) - F$ is ample on $W$ thanks to the fact
that
$-F$ is ample for $\nu $. In order to measure  the
positivity
$\JJ$ with respect to $H$, we ask how small one can
take $s$ to be while keeping the class in question
non-negative: 
\begin{definition} \label{s-inv} The
\textit{$s$-invariant}  of
$\JJ$ with respect to $H$ is defined to be the
positive real number
\[  s_H(\JJ) \ = \ \min \big \{ s \in \RR \ | \ 
\nu ^*( s   H ) - F
\ \text{is nef} \ \big \}.
\] Here $ \nu^* (s  H) - F$ is considered as an
$\mathbf{R}$-divisor on $W$,\footnote{By an
$\RR$-divisor on a variety $V$ we understand an element
of $\text{Div}(V) \otimes \RR$,
$\text{Div}(V)$ denoting the group of Cartier divisors
on $V$.} and to say that it is nef means by definition
that 
\[ s \cdot \big( \nu^* H \cdot C^\prime \big) \ge
\big( F \cdot C^\prime \big) \] for every effective
curve $C^\prime \subset W$. 
\end{definition}

\begin{remark} \label{Compute.after.pullback}
Suppose that $f : Y \lra X$ is a surjective morphism
of projective varieties  with the property
that $\JJ \cdot \OO_Y =
\OO_Y(-E)$ for some effective Cartier divisor $E$ on
$Y$. Then $f$ factors through $\nu$. Recalling that 
nefness can be tested after pull-back by a
surjective morphism, it follows that
\[ s_H (\JJ) \ = \  \min \big \{ s \in \RR \ | \ 
 f^*( s   H ) - E
\ \text{is nef} \ \big \}.  \qed  \]
\end{remark}

\begin{remark} \textbf{(Seshadri Constants).} \label{S.Const.Rem}
Following
\cite{Demailly} and
\cite{Paoletti.JDG} it would be natural to define the
\textit{Seshadri constant} $\epsilon_H(\JJ)$ of $\JJ$
with respect to $H$ to be the reciprocal
\[ \epsilon_H(\JJ) \ = \ \frac{1}{s_H(\JJ)} .
\] However  Definition \ref{s-inv} is more convenient
for our purposes, and we use a different name in order
to avoid the possibility of confusion. When $\JJ$ is
the ideal sheaf of a point $x \in X$ and $L =
\OO_X(H)$, the invariant $\epsilon(L, x) =
\epsilon_H(\JJ)$ was introduced by Demailly   as
a measure of the local positivity  of $L$ at $x$. The
behavior of these Seshadri constants in this case is
very interesting and they have been the focus of
considerable study (cf. \cite{EL1}, \cite{EKL},
\cite{L}, \cite{Bauer1}, \cite{Bauer2}). When
$\JJ$ is the ideal of a smooth subvariety  the
Seshadri constants
$\epsilon_H(\JJ)$ were studied in the interesting
papers (\cite{Paoletti.JDG},
\cite{Paoletti2}), of Paoletti, who considers
especially  smooth curves in threefolds.  Several of
the results in the present note are simple
generaliztions of statements appearing in and
suggested by Paoletti's work, particularly
\cite{Paoletti.JDG}, \S3.  In the past, however, it
was unclear how to use geometric methods to study
these invariants for arbitrary ideals. One of our main
technical observations is that so long as one is
content with asymptotic statements for powers, one
doesn't need restrictions on the geometry of
$\JJ$. This also motivates our study in the present
section of the algebraic properties of the
$s$-invariant. \qed
\end{remark}

We start by comparing this invariant to the twists
needed to generate
$\JJ$.  As customary, set 
\begin{equation} d(\JJ) \ = \ d_H(\JJ) \ = \ \min
\big \{ \ d \in
\ZZ \ | \ \JJ(dH)
\text{ is globally generated}\
\big
\} .\end{equation}
The following is due to Paoletti:
\begin{lemma}   One has
the inequality
\begin{equation} s_H(\JJ) \ \le \ d_H(\JJ). \notag
\end{equation} More  generally, 
\[  s_H(\JJ) \ \le \ \frac{ d_H(\JJ^p) }{p}  \]
for every integer $p \ge 1$. \end{lemma}
\begin{proof} In fact, supppose that $\JJ(dH)$ is
globally generated. Then
\[ \nu^{-1}\JJ (dH)     \ = \
\OO_{V}\big(\nu^*(d H) - E \big) \] is likewise
globally generated and hence nef. Therefore $s_H(\JJ)
\le d$. The second assertion is proven similarly. 
\end{proof}
\begin{remark} We will see later (Theorem
\ref{asympt.reg.thm}) that if
$H$ is very ample, then in fact
$s_H(\JJ) = \lim 
\frac{d_H(\JJ^p)}{p}$. \qed
\end{remark}
\begin{example} \textbf{(Schemes cut out by
quadrics).} Take
$X =
\PP^n$ and
$H$  a hyperplane, and suppose that
$\JJ$ is generated by quadrics, i.e. that $\JJ(2)$ is
spanned by its global sections. Assume in addition
that the zero-locus $Z =
\text{Zeroes}(\JJ)$ is  not  a linear
space. Then
$s_H(\JJ) = 2$. Indeed, the previous Lemma shows that
$s_H(\JJ) \le 2$, and by taking $C^\prime$ in
(\ref{s-inv}) to be the proper transform of  a
general secant line to $Z$, one sees that $s_H(\JJ)
\ge 2$.  
\end{example}

\begin{example} \textbf{(Irrational
$\mathbf{s}$-invariants).} \label{D.C.Const} A
construction used on several occasions by the first
author (cf. \cite{Cut.ZD}) leads to examples where $s_H(\JJ)$ is
irrational. This of course also gives examples where
$s_H(\JJ) < d_H(\JJ)$.\footnote{See also Example
\ref{Pathology.Example}.}  Take
$A$ to be an abelian surface with Picard number
$\rho(A) \ge 3$ (for example $A$ might be the product
of two copies of an elliptic curve). Denote by 
$\text{Nef}(A) \subset NS(A)_{\RR}$ the cone of
numerically effective real divisor classes. Then, as
on any abelian surface, 
\[  \text{Nef}(A) \ = \ \{ \alpha \in NS(A)_{\RR} \ |
\ ( \alpha^2) \ge 0 \ , \ (\alpha
\cdot h ) \ge 0 \ \}, \]
$h$ being any ample class. But the Hodge Index theorem
shows that the intersection form has type $(+ , - ,
\ldots, -)$ on $NS(A)_{\RR}$, and therefore
$\text{Nef}(A)$ is a circular cone.   At least on
suitable $A$,  can then find an effective curve
$C \subset A$, plus an ample divisor class
$H$ such that the ray passing through $-[C]$ in the
direction of $[H]$ meets the boundary of
$\text{Nef}(A)$ at an irrational point, i.e. \[
\inf\big\{s > 0 \ | \ sH - C \in \text{Nef}(A) \big \}
\not \in \QQ . \] Taking $\JJ = \OO_A(-C)$, this means
that $s_H(\JJ)$  is irrational. Note that one can
replace $H$ by $aH$ and $C$ by $C + bH$ ($a, b 
\in \mathbf{N}$), and so arrive at examples with
$C$ and $H$ arbitrarily positive. \qed
\end{example}

\begin{remark}\textbf{(Paoletti's geometric
interpretation of the $\mathbf{s}$-invariant).} 
\label{Paoletti.Interpretation} Suppose that
$Y
\subset X$ is a smooth subvariety with
normal bundle
$N = N_{Y/X}$, and set $\JJ =
\II_{Y/X}$. Then Paoletti \cite{Paoletti.JDG}, p.
487,  shows that the
$s$-invariant
$t_H(\JJ)$ has a simple geometric interpretation, as
follows. Consider first  a non-constant mapping $f :
C \lra X$ from a smooth curve to $X$. If $f(C) \not
\subseteq Y$, then
$f^{-1} \JJ \subset \OO_C$ is an ideal of finite
colength in $\OO_C$, and we define
\[ s^\prime_H(\JJ) \ = \ \underset{\substack{f : C
\lra X \\ f(C) \not
\subseteq Y}}{\text{sup}}
\Big  \{
\frac{\text{colength}( f^{-1}  \JJ )}{ ( C \cdot_f H 
) } \  \Big \} , \] where $(C \cdot_f H)$ denotes the
degree of the divisor $f^* H$ on $C$. Next, put
\[ s^{\prime \prime}_H(\JJ) \ = \ \text{inf} \big \{ \
s > 0 \ | \ N^*(sH) \ \text{is nef}
\ \big \}, \] the nefness of a bundle twisted by an
$\QQ$ or $\RR$ divisor being defined in the evident
manner (cf. \cite{PAG}, Chapter 2). Then
\[ s_H(\JJ)  \ = \ \max \, \{ \ s^\prime_H( \JJ) \ , \
s^{\prime
\prime}_H(\JJ) \ \}. \] In fact, given $f : C \lra X$
as above, let
$f^\prime : C \lra W$ be the proper transform of
$f$. Then
$\text{colength}(f^{-1} \JJ ) = \big( C
\cdot_{f^\prime} F \big) $, and consequently
$s_H^{\prime}(\JJ)$ is the least real number $s^\prime
> 0$ such that
$\nu_0^*(s^\prime H) - F$ has non-negative degree on
every  curve
$C^\prime \subset W$ not lying in the exceptional
divisor $F \subset W$. Similarly,
$s^{\prime \prime}_H(\JJ)$ controls the nefness of
$\OO_{F}\big(
\nu_0^*(s^{\prime \prime} H) - F\big)$.

\end{remark}

 For constructing examples, it is useful to
understand something about how the $s$-invariant
behaves in ``chains" of subvarieties. With
$X$ as before, consider then a sequence non-singular
irreducible subvarieties. 
\[ Z \subseteq Y \subseteq X,   \]  and fix an ample
divisor $H$ on $X$. There are three naturally defined
ideal sheaves in this setting, and we can consider
the correponding
$s$-invariants
\[ s_H(\II_{Z/X}) \ \ , \ \ s_H(\II_{Y/X})  \ \
\text{and} \ \ s_H(\II_{Z/Y}); \] in the third case we
view $H$ as an ample divisor on
$Y$, and compute on the blow-up of 
$Y$.  One evidently has the inequality $s_H(\II_{Z/Y})
\le s_H(\II_{Z/X})$, and in favorable situations the 
two invariants in question coincide: 
\begin{proposition} \label{Chains} In the situation
just described, assume that 
$s_H(\II_{Y/X})
< s_H(\II_{Z/Y})$. Then 
\[ s_H(\II_{Z/X}) = s_H( \II_{Z/Y}). \] 
\end{proposition} 
\begin{proof} We keep the notation introduced in
Remark \ref{Paoletti.Interpretation}. It is evident
that
\[
s_H^\pr(\II_{Z/X}) \   \le \  \max \{ \,
s_H^\pr(\II_{Z/Y}) \, , \, s_H^\pr (\II_{Y/X}) \, \},
\]
and it follows from the
conormal bundle sequence $ 0 \rightarrow N^*_{Y/X} | Z
\rightarrow N^*_{Z/X} \rightarrow N^*_{Z/Y}
\rightarrow 0$ that
\[ s_H^{\pr\pr}(\II_{Z/X})
\  \le \ \max
\{ \, s_H^{\pr \pr}(\II_{Z/Y}) \, , \, s_H^{\pr\pr}
(\II_{Y/X}) \, \} . \]
The assertion is then a consequence of Remark
 \ref{Paoletti.Interpretation}.
 \end{proof}

\begin{example} \textbf{(Irrational
$\mathbf{s}$-invariants on projective space).}
One can combine Example
\ref{D.C.Const} with the previous Proposition
 to arrive at a quick  example of a curve
$C \subset \PP^r$, with ideal sheaf $\JJ =
\II_{C/\PP^r}$, such that $s_H(\JJ)$ is irrational,
$H$ being the hyperplane class.  Specifically, take a
very ample divisor $H$ on an abelian surface $A$, plus
a curve $C \subset X$, such that
$s_H(\OO_A(-C))$ is an irrational number $ > 2$, and
such that $A$ is cut out by quadrics under the
embedding $A
\subset \PP = \PP^r$ defined by $H$.  [Starting with
any $C$ and
$H$ giving irrational invariant, first replace $H$ by
$4H$ to ensure that $\II_{A/\PP}$ is generated by
quadrics, and then replace $C$ by $C + mH$ for $m \gg
0$ to make
$s_H(\OO_A(-C)) > 2$.] Then
$s_H(\II_{A/\PP}) = 2$, so by applying the previous
example to the chain $C
\subset  A \subset
\PP^r$, we find that  $s_H(\JJ) = s_H(\OO_A(-C))$.
Examples of curves in $\PP^3$ having irrational
$s$-invariant were given by the first author in
\cite{Cut.Irr}\footnote{The cited paper deals with curves
having irrational asymptotic
Castelnuovo-Mumford regularity, but Theorem
\ref{asympt.reg.thm} shows that this is the same as
irrational
$s$-invariant.}, but they involved more computation.
\qed
\end{example}

We conclude this section with a result giving some
algebraic properties of the
$s$-invariant:
\begin{proposition} Let $X$ be an irreducible
projective variety, and
$H$ an ample divisor on $X$.
\begin{enumerate}
\item[(i).]  Given ideal sheaves $\JJ_1,\JJ_2 
\subseteq \OO_X$, one has the inequalities:
\begin{align*} s_H\big( \JJ_1 \cdot \JJ_2 \big) \ &\le
\ s_H(\JJ_1) + s_H(\JJ_2) \\
s_H\big( \JJ_1 + \JJ_2 \big) \ & \le
\ \max \, \big \{ \, s_H(\JJ_1) \, , \,  s_H(\JJ_2)
\, \big \}.
\end{align*}
\item[(ii).]  If ${\bar \JJ} \subseteq \OO_X$ denotes
the integral closure of $\JJ$, then 
\[ s_H({ \bar  \JJ}) = s_H(\JJ). \] 
\end{enumerate}
\end{proposition}
\noi For basic facts about the integral closure of an
ideal, see
\cite{Teissier}.
\begin{proof} We will apply Remark
\ref{Compute.after.pullback}. Thus for (i), let
$f : Y \lra X$ be a surjective mapping from an
irreducible  variety
$Y$ which dominates the blowings-up of $X$ along
$\JJ_1$,
$\JJ_2$ and $\JJ_1 + \JJ_2$. Thus $Y$ carries
effective Cartier divisors $E_1, E_2$ and $E_{12}$
characterized by
\[ \JJ_1 \cdot \OO_Y = \OO_Y(-E_1) \ \ , \ \ \JJ_1
\cdot
\OO_Y = \OO_Y(-E_1) \ \ , \ \ \big(\JJ_1+\JJ_2 \big)
\cdot \OO_Y =
\OO_Y(-E_{12}). \]
Note that then 
\begin{equation}
\big( \JJ_1  \JJ_2 \big) \cdot
\OO_Y = \OO_Y\big( - (E_1 + E_2) \ \big). \tag{*}
\end{equation}

Write $s_1 = s_H(\JJ_1)$ and $s_2 = s_H(\JJ_2)$. Then
$f^*(s_1H) - E_1$ and $f^*(s_2H) - E_2$ are nef on
$Y$, and consequently so is their sum $f^*\big( (s_1
+ s_2)H\big) - (E_1 + E_2) $. The first inequality in
(i) then follows from (*). For the second, set $s =
\max
\{ s_1, s_2 \}$ and note that one has a surjective
map
\[ \OO_Y( - E_1) \oplus \OO_Y( - E_2) \lra \OO_Y(
-E_{12}) \]
of vector bundles on $Y$. By definition of $s$, the
bundle on the left becomes nef when twisted by the
$\RR$-divisor
$f^*(sH)$. Since quotients of nef bundles are nef, 
this implies that 
$f^*(sH) - E_{12}$ is  nef, and the required
inequality follows.\footnote{We are implicitly using
here the fact that nefness  makes sense for twists of
bundles by $\QQ$- or
$\RR$-divisors, and that the usual formal properties
are satisfied. These facts are worked out in Chapter
2 of the forthcoming book \cite{PAG}, but the reader
can easily verify the required assertion directly by
considering the evident $\RR$-divisors on the
projectivization $\PP\big( \OO_Y( - E_1) \oplus
\OO_Y( - E_2) \big) \lra X$. }

For (ii), we use the fact (cf. \cite{Teissier}, p.330)
that 
$\Bl_{\JJ}(X)$ and 
$\Bl_{\bar \JJ}(X)$ have the same normalization
$V$, which sits in a commutative diagram:
\[
\xymatrix{
& V \ar[dl]  \ar[dd]^{\mu} \ar[dr]  \\
\Bl_{\bar \JJ}(X) \ar[dr]_{\bar \nu } & & \Bl_{\JJ}(X)
\ar[dl]^{\nu }
\\ & X 
}
\]
 Moreover, the exceptional divisors $F$ and ${\bar
F}$ of $\nu $ and ${  \nu }$ pull back to the
same divisor $E$ on $V$. Invoking again Remark
\ref{Compute.after.pullback} one has
\[ s_H(\JJ) \ = \ \inf \{ s > 0 \ | \ \mu^*(sH) - E \
\text{is nef } \} \
 = \ s_H(\bar \JJ), \] as required. \end{proof}

\vskip 10pt

\section{Degree and Nilpotency Bounds}

In the present section, we show how the $s$-invariant
governs bounds on the degrees of zeroes of an ideal
and its index of nilpotency. For the most part this
involves only small modifications to computations
appearing for instance in \cite{EL99}, so we shall be
brief.

We start by fixing some additional notation. Let $X$
be a non-singular  irreducible quasi-projective
variety of dimension
$n$ --- which for the moment we do not assume to be
projective --- and suppose that $\JJ
\subset
\OO_X$ is a coherent sheaf of ideals on $X$. As
before we denote by 
$\nu : W = \Bl_{\JJ}(X)
\lra X$ the blowing up of $\JJ$, with exceptional
divisor $F$.  Consider now the normalization $p : V
\lra W
$ of
$W$, with
$\mu : V
\lra X$ the natural composition:
\[
\xymatrix{
V \ar[r]_p \ar@/^1pc/[rr]^{\mu} & W \ar[r]_{\nu} &
X.  } \]
We denote by  $E = p^* F$ the pull-back to $V$ of the
exceptional divisor $F$ on $W$. Thus $E$ is an
effective Cartier divisor on $V$, and 
\[ \JJ \cdot \OO_V \ = \ \OO_V(-E) . \]
Note that $\OO_V(-E)$ is ample relative to $\mu$, and
in particular is ample on every fibre of $\mu$.

Now $E$ determines a Weil divisor on $V$, say
\[ [ E ]  \ = \ \sum_{i=1}^{t} \ r_i \cdot [E_i] , \]
where the $E_i$ are the irreducible components of the
support of $E$, and $r_i > 0$. Set
\[ Z_i  \ = \ \mu (E_i) \ \subseteq \ X, \]
so that $Z_i$ is a reduced and irreducible subvariety
of $X$. Following \cite{Fulton}, the $Z_i$ are called
the \textit{distinguished subvarieties} of $\JJ$.
(Note that several of the $E_i$ may have the same
image in $X$, in which case there will be repetitions
among the $Z_i$.  However this doesn't cause any
problems.)  Denoting by 
\[ Z \ = \ \text{Zeroes} \big( \sqrt{\JJ} \big) \]
the reduced zero-locus of $\JJ$, one has then the
decomposition
\[ Z   \ = \ \cup Z_i \]
of $Z$ as a union of distinguished subvarieties. Thus
each irreducible component of $Z$ is distinguished, but
there can be ``embedded" distinguished subvarieties 
as well. We refer to the positive
integer $r_i$ as the coefficient attached to
$Z_i$, and we define
\begin{equation} r(\JJ) \ =_{\text{def}}   \ \max
\{ r_i \} . \end{equation}

The following result, implicit in \cite{EL99} and
independently observed by Hickle \cite{Hickle}, shows that
the invariant
$r(\JJ)$ controls the index of nilpotency of $\JJ$: 
\begin{theorem} \label{Abstract.NSS}  One has
\[ \Big ( \sqrt{\JJ} \Big ) ^{n \cdot r(\JJ)}
\subseteq
\JJ.\] More generally, $\big( \sqrt{\JJ}
\big)^{(n+1-p)\cdot r(\JJ)} \subseteq \JJ^p$ for every
integer
$p
\ge 1$. 
\end{theorem}
\begin{proof}[Sketch of Proof] 
One checks right away as in \cite{EL99}, (2.1) and
(2.4), that
\[
\big(\sqrt{\JJ}\big)^{\ell \cdot r(\JJ)}
\ \subseteq  \ \mu_* \OO_V(-\ell E) \ = \ \overline{
\JJ^\ell}\] 
for every $\ell \ge 0$. The stated inclusions then
follow from the Brian\c con--Skoda theorem
(cf. \cite{Huneke}).
\end{proof}

In order to give Theorem \ref{Abstract.NSS} some real 
content, one needs an upper bound on $r(\JJ)$. It
would be interesting to know whether one can give
useful statements in a purely local setting. However
globally they follow (Corollary
\ref{Bezout.1}) from the  fact 
 that one has  Bezout-type inequalities for the
degrees of the distinguished subvarieties in terms of
the
$s$-invariant of $\JJ$.

Assume henceforth that $X$ is projective, and fix an
ample divior class $H$ on $X$. 
\begin{proposition} \label{Deg.Bound} Let $s =
s_H(\JJ)$ be the $s$-invariant of $\JJ$ with respect to
 $H$.  Then 
\[
\sum_{i = 1}^{t}  \ r_i  \cdot s^{\dim Z_i} \cdot
\deg_H  Z_i \ \le \ s^n \cdot \deg_H X, \]
where for any subvariety $V \subseteq X$, $\deg_H V =
\big( H^{\dim V} \cdot V \big)$ denotes the degree of
$V$ with respect to $H$. 
\end{proposition}

\begin{corollary}\label{Bezout.1} In the situation of
the Proposition one has
 \[ \sum s^{\dim Z_i} \deg_H
Z_i \ \le \ s^n \ \deg_H X, \]
 and the integer
$r(\JJ) = \max \{ r_i \}$  satisfies \[ r(\JJ) \le
s_+^n
\cdot \deg_H X , \]
where $s_+ = \max \{ 1, s_H(\JJ) \}$. \qed
\end{corollary}

The Proposition can be deduced from general positivity
results  due to Fulton and the
third author \cite{FL}.  However  following \cite{EL99}
we indicate a direct proof using classical
intersection theory.
\begin{proof}[Sketch of proof of Proposition
\ref{Deg.Bound}] Consider the classes
\[ h  \ = \ [ \mu^* H ] \ \ , \ \ m \ = \ [\mu^*(sH) -
E] \ \in \text{NS}(V)_{\RR} \]
in the vector space of numerical equivalence classes
on $X$ with real coefficients.  Thus $m$ is a nef
class --- so in particular $\int_V m^n \ge 0$ and
$\int_{E_i} (s \cdot h)^j \cdot m^{n-1-j} \ge 0$ for
all $i$ and $j$ --- and
$[E] = s
\cdot h - m$.  Arguing as in the proof of
Proposition 3.1 in
\cite{EL99}, one then finds that
\[
\begin{aligned}
 s^n \cdot \deg_H (X) &= \int_V \big(s\cdot h\big) ^n
\\
 &\ge \int_V \big (  ( s \cdot h  ) ^n - m^n \big
)
\\
 &= \int_V \Big( (s \cdot h) - m \Big) \Big(  \sum_{j
= 0}^{n-1} \ (s \cdot h)^{j}\cdot m^{n-1-j} \Big) \\
&= \int_{[E]} \Big(  \sum_{j
= 0}^{n-1} \ (s \cdot h)^{j} \cdot m^{n -1-j} \Big) 
\\  &\ge \sum_{i = 1}^t r_i \cdot \int_{E_i}
(s\cdot h)^{\dim(Z_i)} \cdot m^{n-1-\dim(Z_i)}, \\
&\ge \sum_{i = 1}^t r_i \cdot s^{\dim Z_i} \cdot
\deg_H Z_i,
\end{aligned}
\]
as required.
\end{proof}

\begin{remark}
Suppose that $Z_1, \ldots , Z_p$ are the (distinct)
irreducible components of $Z$.  Then arguing as in
\cite{Fulton} (4.3.4) and (12.2.9) one finds
\[ \sum_{i = 1}^p \,  e_{Z_i} (\JJ) \cdot s^{\dim Z_i}
\cdot \deg_H Z_i \ \le \  s^n \cdot \deg_H X, \]
where $e_{Z_i}(\JJ)$ is the Samuel multiplicity of
$\JJ$ along $Z_i$. \qed
\end{remark}

One does not expect Bezout-type inequalities such as
\ref{Bezout.1} to capture the embedded components of
$\JJ$ in the sense of primary decomposition (see
\cite{BM}, \cite{Flenner}, \cite{Kollar} and
Example \ref{Pathology.Example} below). Somewhat unexpectedly,
however,  the situation is different when
$\JJ$ is integrally closed:
\begin{corollary} Assume that $\JJ$ is integrally
closed, and let $Y_1, \ldots , Y_q \subseteq X$ be 
the irreducible subvarieties defined by all the
associated primes of $\JJ$ (minimal or embedded).
Then 
\[ \sum_{j = 1}^q  s^{\dim Y_j} \deg_H
Y_j \ \le \ s^n \ \deg X.\]
\end{corollary}

\begin{proof} It is enough to show that every
associated subvariety is distinguished. To this end,
let
$\frq_i = \mu_*
\OO_Y\big( - r_i   E_i  \big)$ be the sheaf of all
functions on $X$ whose pull-backs to $V$  vanish
to order
$\ge r_i$ along the Weil divisor $E_i$. Then $\frq_i
\subseteq
\OO_X$ is a primary ideal, and one has 
\begin{align*} \overline \JJ \ &= \ \mu_* \OO_Y(- E)
\\ &= \  \bigcap_{i=1}^t \, \mu_* \OO_Y \big(- r_i
 E_i
\big).
\end{align*}
Since $\JJ$ is integrally closed this means that we
have the (possibly redundant) primary decomposition
$\JJ = \cap
\frq_i$. In particular every
associated prime of
$\JJ$ must occur as the radical of one of the
$\frq_i$, i.e. as one of the distinguished
subvarieties. 
\end{proof}

\begin{remark}
The argument just given to show that each associated
subvariety of $\overline \JJ$ is distinguished
appears a number of times in the literature (e.g.
\cite{Hickle}). However it seems to have been
overlooked that this leads to degree bounds on
associated subvarieties for integrally closed ideals.
\qed
\end{remark}

\begin{remark} The same argument shows more generally
that for any ideal $\JJ$, the bound $\sum 
s^{\dim Y} \deg_H Y \le   s^n \ \deg_H X$ holds if
one sums over all subvarieties $Y$ defined by an
associated prime ideal of the integral closure
$\overline {\JJ^p}$ of some power of 
$\JJ$. \qed
\end{remark}

\begin{example} \textbf{(Pathological ideals with
fixed $\mathbf{s}$-invariant).}
\label{Pathology.Example} We construct here a  family
of ideals having fixed
$s$-invariant but arbitrarily many embedded points.
The same examples will show that the regularity bounds
presented in the next section only hold
asymptotically. For simplicity we work over the
complex numbers $\CC$, but in fact one could deal
with an arbitrary algebraically closed ground-field.

In order to highlight the underlying geometric
picture, we start with a local discussion. Working in
affine three-space $X = \AAA^3$ with coordinates $x,
y, t$, consider the ideal $\fra = \fra_p = ( x^2,
p(t)\cdot xy, y^2)
\in
\CC[x,y,t]$, where $p(t) \in \CC[t]$ is a polynomial
in $t$.  Then the zeroes of $p(t)$ along the line
$\ell$ defined by $ \{ x = y = 0 \}$ are embedded
points of
$\fra$. On the other hand, let
$f: Y =
\Bl_{(x,y)}(X)
\lra X$ is the blowing up of $X$ along $\ell$, with
exceptional divisor
$E$. Then one checks that $\fra \cdot \OO_Y =
\OO_Y(-2E)$, so in other words on $Y$ the ideal
$\fra$ cannot be distinguished  from the square
$(x,y)^2$ of the ideal
of $\ell$.\footnote{Geometrically, the important point
is that for every complex number
$a
\in \CC$, the homogeneous polynomials
$x^2 \, ,\, p(a)xy\, , \, y^2 \in \Gamma \big (\PP^1,
\OO_{\PP^1}(2) \big)$ span a base-point linear
series. More algebraically,
observe that already $(x^2, y^2) \cdot \OO_Y =
\OO_X(-2E)$.} The idea is that in the global setting,
the $s$-invariant will be  computed on the blow-up of
the line (Example \ref{Compute.after.pullback}), and
so cannot detect the embedded points.

This local construction is easily globalized. Take $X
= \PP^3$ with homogeneous coordinates $X, Y, Z, W$,
 fix a homogeneous polynomial $P_d = P_d(Z,W) \in
\CC[Z,W]$ of degree $d$, and let $\JJ = \JJ_P
\subseteq
\OO_{\PP^3}$ be the ideal sheaf spanned by the
homogeous polynomials $X^2$,
$P_d \cdot XY$ and $Y^2$. Denoting by $L \subseteq
\PP^3$ the line $\{ X = Y = 0 \}$, one sees as above
that the zeroes of $P_d$ along $L$ are embedded
points of $\JJ$, so for general $P$ there will be $d$
such.  As before let $Y = \Bl_L(\PP^3) \lra \PP^3$ be
the blowing-up of $L$, with exceptional divisor $E$. 
Then $\JJ_P \cdot \OO_Y = \OO_Y( -2E )$, so it
follows from Example \ref{Compute.after.pullback}
that $s_H(\JJ_P) = 2$ for every $P$ ($H$ being the
hyperplane divisor). In particular, the number of
embedded points cannot be bounded in terms of the
$s$-invariant. \qed
\end{example}

\section{Asymptotic Regularity and Degree Bounds}

In the present section we bound the
"complexity" of large powers of an ideal sheaf in
terms of its $s$-invariant. 

As above, let $X$ be a non-singular irreducible
projective variety of dimension $n$"  We assume in
this section that
$H$ is a \textit{very} ample divisor on $X$. In this
case the   Castelnuovo-Mumford regularity of a
coherent sheaf $\FF$ on $X$ is defined just as in the
classical setting of projective space:
\begin{definition} A coherent sheaf $\FF$ is
$m$-regular (with respect to $H$) if
\[ \HH{i}{X}{\FF\big((m-i)H\big)} = 0 \for i > 0.
 \]
The \textit{regularity} $\reg_H(\FF)$ of $\FF$ is the
least integer $m$ for which $\FF$ is
$m$-regular.\footnote{If $\FF$ is $m$-regular for
every $m \in \ZZ \ $ --- which will occur if and only
if $\FF$ is supported on a finite set --- we put
$\reg_H(\FF) = - \infty.$}
\qed
\end{definition}
\noi Just as in the classical case, if $\FF$ is
$m$-regular for some integer $m$, then $\FF(mH)$ is
globally generated, and
$\FF$ is also
$(m + 1)$-regular. We view the regularity of a sheaf
as a measure of its cohomological complexity.  When
$X = \PP^n$, this regularity has a well-known
interpretation as bounding the degrees of the
generators of the modules of syzygies of the
 module corresponding to $\FF$ (see \cite{BM}). 

Fix now an ideal sheaf $\JJ \subset \OO_X$ with
$s$-invariant $s_H(\JJ)$.  As above we denote by
$d_H(\JJ^p)$ the least integer $d \ge 0$ such that
$\JJ^p(dH)$ is globally generated.
\begin{theorem} \label{asympt.reg.thm} The quantities
$\frac{\reg_H(\JJ^p)}{p}$ and
$
\frac{d_H(\JJ^p)}{p}$ tend to   limits as $p \to
\infty$, and one has:
\[ \lim_{p \to \infty} \frac{ \reg_H(\JJ^p)}{p}   \ =
\lim_{p \to \infty} \frac{ d_H(\JJ^p)}{p} \ = \ 
s_H(\JJ).
\]
\end{theorem}
\begin{proof}
Set $d_p = d_H(\JJ^p)$ and $r_p = \reg_H(\JJ^p)$. Note
to begin with that $d_{\ell + m} \le d_\ell + d_m$
for all $\ell , m \ge 0$, from which it follows that
the limit $\lim_{p \to \infty} \frac{d_p}{p}$ exists.
  Call this limit $\bar d$. We will prove the
theorem by establishing (from right to left) the
inequalities
\begin{equation} \label{reg.eqn}
\limsup \frac{r_p}{p} \ \le \ s_H(\JJ) \ \le \ \bar d
\ \le \ \liminf \frac{r_p}{p}. 
\end{equation} 

Starting with the right-most inequality in 
(\ref{reg.eqn}), recall that if $\JJ^p$ is
$m$-regular with respect to $H$ then $\JJ^p(mH)$ is
globally generated.  Therefore
$d_p \le r_p$ for every $p > 0$, and in particular
$\lim \tfrac{d_p}{p} \le
\liminf \frac{r_p}{p}$.

We next show that $s_H(\JJ) \le \lim \frac{d_p}{p}
= \bar d$. To this end, fix any $\epsilon > 0$. Then
we can choose large positive integers $p_0 , q_0 > 0$
such that \[\frac{d_{p_0}}{p_0} \ \le \
\frac{q_0}{p_0}
\ \le \ \bar d +
\epsilon, \] so that in particular $\JJ^{p_0}(q_0H)$
is globally generated. Writing as before $\nu : W =
\Bl_{\JJ}(X) \to X$ for the blow-up of $\JJ$, with
exceptional divisor $F$, it follows that $\nu^*(q_0H)
- p_0F$ is globally generated and hence nef. Therefore
$s_H(\JJ) \le \tfrac{q_0}{p_0} \le \bar d +
 \epsilon$,
as required. 

It remains to prove that $\limsup
\frac{r_p}{p} \le s_H(\JJ)$. To this end we use a
theorem of Fujita \cite{Fujita} to the effect that
Serre Vanishing remains valid even after twisting by 
arbitrary nef divisors. Specifically, consider an
irreducible projective variety $V$, and fix an ample
divisor $A$ plus a coherent sheaf $\FF$ on $V$.
Fujita shows that there is an integer $m_0 = m_0(A,
\FF)$ such that for any nef divisor $B$:
\begin{equation} \HH{i}{V}{\FF(mA + B)} = 0 \fall i >
0 \ \text{and} \ m
\ge m_0. \tag{*} \end{equation}
(The important point here is that $m_0$ is
independent of $B$.)

We propose to apply (*) on the blowing-up $W
= \Bl_{\JJ}(X)$ of
$\JJ$. Given $\epsilon > 0$, choose large
integers $q_0, p_0$ such that
\[ s_H(\JJ) \ < \ \frac{q_0}{p_0} \ < \ s_H(\JJ) +
\frac{\epsilon}{2} . \]
Then $\nu^*(q_0H) - p_0F$ is ample, so there exists
an integer $m_0$ such that if $m \ge m_0$ then for any
nef divisor $P$ on $W$, the bundles associated to the
divisors $\nu^*(mq_0H) - mp_0F + P$ have vanishing
higher cohomology. Now fix any integer $p \ge m_0
p_0$, and write 
\[ p = m p_0 + p_1 \ \ \text{with} \ \ 0 \le p_1 <
p_0 \ \text{and} \  m \ge m_0.\] 
Then $\nu^*(q_0H) -
p_1F$ is nef (in fact ample), and consequently we
have the vanishing of the higher cohomology of the
line bundle
\[  \OO_W \big(  \, \nu^*\left ( ( m + 1)q_0 H \,
\right ) - pF
\, \big ) .  \] It now follows from Lemma
\ref{powers.as.pushforwards} below --- and this is
the crucial point --- that
\[
\HHH{i}{X}{\JJ^p  \big( (m+1)q_0H\big) } = 0 \for i >
0\]  provided that 
$p$ is sufficiently large. Therefore
$\JJ^p$ is $\big( (m+1)q_0 + n \big)$-regular for $p
\gg 0$, and  consequently
\[ \frac{r_p}{p} \ \le \ \frac{ (m+1)q_0 + n}{p} \
\le \ \frac{q_0}{p_0} + \frac{q_0 + n}{mp_0}. \]
By taking $p$ (and hence also $m$) to be large
enough, we can arrange that the second term on the
right is $\le \tfrac{\epsilon}{2}$, so that
$\frac{r_p}{p} \le s_H(\JJ) + \epsilon$ for $p \gg
0$. Therefore $\limsup \tfrac{r_p}{p} \le s_H(\JJ)$,
and we are done. 
\end{proof}

The following Lemma played an essential role in the
proof just completed. It shows that one can realize
large powers of an ideal $\JJ \subseteq \OO_X$
geometrically from the natural divisor on the
blow-up.\footnote{If $\JJ$ defined a smooth
subvariety of $X$, then the corresponding statement
would be true for all powers.} This fact is surely
not new, but we include a proof for the convenience
of the reader. 
\begin{lemma}  \label{powers.as.pushforwards} Let  
$\JJ \subseteq \OO_X$ be an ideal sheaf on
$X$, and \[
\nu : W = \Bl_{\JJ}(X) \lra X \] the
blowing-up of $\JJ$, with exceptional divisor $F$.
There exists an integer $p_0 > 0$ with the property
that if $p \ge p_0$, then
\begin{equation}\nu_* \OO_W \big( -pF \big) =
\JJ^p , \tag{*} \end{equation}
 and for any divisor $D$ on
$X$:
\[ \HH{i}{X}{\JJ^p(D)}
= \HH{i}{W}{\OO_W(\nu^*D - pF)}  \] 
for all $i \ge 0$. 
\end{lemma}

\begin{proof} Since $\OO_W(-F)$ is ample for $\nu$,
it follows from Grothendieck-Serre vanishing that
 \[ R^j \nu_* \OO_W(-pF) = 0 \
\for  j > 0 \ \text{and} \ p \gg 0.\] The isomorphism
on global cohomology groups is then a consequence of
(*) thanks to the Leray spectral sequence. 

As for (*),
the assertion is local on $X$, so we may assume that
$X$ is affine. Choosing generators $g_1 , \ldots ,
g_r \in \JJ$ gives rise to a surjection $\OO_X^r \lra
\JJ$, which in turn determines an embedding
\[  W = \Bl_{\JJ}(X) \ \subseteq \ \PP \big( \OO_X^r
\big) = \PP^{r-1}_X \]
in such a way that $\OO_{\PP^{r-1}_X}(1) \mid W =
\OO_W(-F)$. Write $\pi : \PP^{r-1}_X = X \times
\PP^{r-1} \lra X$ for the projection. Serre vanishing
for
$\pi$, applied  to the ideal sheaf
$\II_{W/\PP^{r-1}_X}$, shows that if
$p \gg 0$  then the natural  homomorphism
\begin{equation}\pi_* \OO_{\PP^{r-1}_X}(p) \lra \pi_*
\OO_W(-pF)
\tag{**} \end{equation} is surjective. On the other
hand, recalling that $\pi_* \OO_{\PP^{r-1}_X}(k) =
S^k \big( \OO_X^r \big)$ for every $k \ge 0$, one sees
that the image of (**) is exactly $\JJ^p$. It follows
that
$\nu_* \OO_W
\big( -pF \big) = \JJ^p$ for $p \gg 0$, as asserted.
\end{proof}

\begin{remark}
The use of Serre Vanishing in the proof of Theorem
\ref{asympt.reg.thm} was suggested by Demailly's
proof of Theorem 6.4 in \cite{Demailly}. Proposition
3.3 of
\cite{Paoletti.JDG} uses a similar argument to prove
a result for zero-loci of vector bundles that is
rather close in spirit to \ref{asympt.reg.thm}. \qed
\end{remark}

Finally, we turn to asymptotic bounds on the
arithmetic degree of $\JJ^p$. In a general way we
follow the approach of Bayer and Mumford, suitably
geometrized. We start by recalling the definition of
the arithmetic degree from the viewpoint of
\cite{Kollar}. 

Assume then that $X$ is a non-singular irreducible
projective variety of dimension $n$ carrying a fixed
ample divisor class $H$, and let $\FF$ be a coherent
sheaf on $X$. Then there is a canonical filtration
\[
0 \subseteq \FF^n \subseteq \FF^{n-1} \subseteq \ldots
\subseteq \FF^{1} \subseteq \FF^0 = \FF \]
where $\FF^k \subseteq \FF$ is the subsheaf
consisting of sections whose support has codimension
$\ge k$ in $X$.  As in \cite{Fulton}, Example
18.3.11, one can in a natural way associate to the
quotient
$\FF^k / \FF^{k+1}$ a codimension $k$ cycle $[\FF^k /
\FF^{k+1}] \in Z^{k}(X)$, and then the codimension
$k$ contribution to the arithmetic degree of $\FF$ is
defined to be
\[ 
\adeg_H^k \big ( \FF \big) \ = \ \deg_H \big( \,
[\FF^k /
\FF^{k+1}] 
\, \big ), \]
where as indicated the degree of $[\FF^k /
\FF^{k+1}]$ is  computed with
respect to the fixed polarization $H$. For an ideal
$\JJ
\subseteq
\OO_X$ one sets $\adeg^k_H \big(\JJ\big) = \adeg^k_H
\big( \OO_X / \JJ\big)$.  Thus $\adeg^k(\JJ)$
measures the degrees of the codimension $k$
components of $\JJ$ (both mimimal and embedded),
counted with  suitable multiplicities. 

A variant of the following Lemma was implicitly
used by Bayer and Mumford in a similar context, and
re-examined in \cite{MV} . 
\begin{lemma} \label{adeg.Bertini} Still assuming
that $H$ is very ample, let
$D \subseteq X$ be a general divisor linearly
equivalent to
$H$, and let $\FF_D =
\FF \otimes_{\OO_X} \OO_D$ denote the restiction of
$\FF$ to $D$. If $k \le n-1$ then 
\[ adeg^k_H \big( \, \FF_D \, \big) \ = \ \adeg^k_H
\big( \, \FF \, \big) ,
\] where the degree on the left is computed with
respect to the ample line bundle $\OO_D(H)$ on $D$.
\end{lemma}

\begin{proof} [Indication of Proof] The essential
point is to show that if $\mathcal{M}$ is an
equidimensional $\OO_X$-module without embedded
components, then the restriction $\mathcal{M}_D$ of
$\mathcal{M}$ to $D$ is also equidimensional without
embedded components (see  \cite{MV} for an argument in
a similar setting). Once one knows this, one can
deduce the lemma from the fact \cite{Fulton},
Examples (18.3.6) and (18.3.11),  that
$\deg_H
\big( \, [\FF^k /
\FF^{k+1}] 
\, \big )$ governs the leading term of the Hilbert
polynomial of the sheaf in question. We leave details
to the reader. \end{proof}

In the spirit of \cite{BM}, Proposition
3.6, we show that --- at least asymptotically --- the
arithmetic degrees of large powers of an ideal are
bounded in terms of their regularity:
\begin{theorem} \label{asympt.adeg.bound}
Suppose as above that $X$ is an irreducible
projective variety, and assume that $H$ is a very
ample divisor on $X$.  Let $\JJ \subseteq \OO_X$ be an
ideal sheaf on
$X$, and set $\overline{\reg}_H(\JJ) = \lim
\tfrac{\reg_H  (\JJ^p)  }{p}$. Then for every  
  $0 \le k \le n$:
\[ \limsup_{p \to \infty} \frac{ \adeg^k_H \big( \JJ^p
\big) }{p^k}
\ \le \ \frac{
\big( \,\overline{reg}_H(\JJ)\big)^k}{k!}
\cdot \deg_H (X).
\]
\end{theorem}

\begin{corollary} In the situation of the
Proposition, 
\[ \limsup_{p \to \infty} \frac{ \adeg_H^k \big( \JJ^p
\big) }{p^k}
\ \le \ \frac{{s_H(\JJ)}^k}{k!} \cdot \deg_H (X). \qed
\]
\end{corollary}

\begin{proof}[Proof of Theorem
\ref{asympt.adeg.bound}]
Let $D \in |H|$ be a general divisor linearly
equivalent to $H$, and consider the restriction
$\JJ_D = \JJ \cdot \OO_D$ of $\JJ$ to $D$. 
According to a theorem of Ratliff \cite{Ratliff}  there
are only finitely many prime ideals which appear as
associated primes for any of the ideals $\JJ^p$ for
$p \ge 1$. So we may assume that $\OO_X(-D)$ does not
contain any of these primes, so that the
sequence 
\[ 0 \lra \JJ^p(-D) \overset{\cdot D}
{\lra} \JJ^p \lra \JJ^p_D \lra 0 \]
is exact for every $p$. This sequence shows that 
$\reg_H(\JJ_D^p) \le \reg_H(\JJ^p)$ for every
$p$, where by abuse of notation we are writing $H$ for
the class of the restriction $ \OO_D(H)$ to $D$.
Consequently
$\overline{\reg}_H(\JJ_D)
\le
\overline{\reg}_H(\JJ)$. Similarly, Lemma
\ref{adeg.Bertini} shows that 
$\adeg_H^k(\JJ^p)
\le
\adeg_H^k(\JJ_D^p)$ for fixed $p$ provided that $k \le
n-1$.  As
we are working over an uncountable ground
field, we can assume by taking $D$ to be very
general that this holds simultaneously for
all $p \ge 1$.\footnote{This is the only
point at which we use the uncountability of
the ground field. One suspects that one could
avoid this.} Since of course also $\deg_H D = \deg_H
X$,   if
$k
\le n-1$ it therefore suffices to prove the
Proposition for
$D$. So by induction on
$n = \dim X$ we can assume that $k = n$.  

Supposing then that $k = n $, 
we need to bound as a function of
$p \gg 0$ the length of the (finitely supported)
subsheaf
\[ \QQQ_p  \ 
\subseteq \ \OO_X / \JJ^p  \]
of sections having zero-dimensional support.
Equivalently, we need to bound for $p \gg 0$ the
dimension $\hh{0}{X}{\QQQ_p}$. 
To this end, observe first of all that for every
integer $q \in \NN$ there is an inclusion
\begin{equation} \label{H0.Qp}
 \HH{0}{X}{\QQQ_p} \ \cong \ \HH{0}{X}{\QQQ_p \otimes
\OO_X(qH)} \ \subseteq \HH{0}{X}{\big( \OO_X / \JJ^p
\big )  \otimes \OO_X(qH)}. 
\end{equation}
The plan is to estimate the dimension of the group on
the right for a suitable  integer $q$. Fix
$\epsilon > 0$ plus large integers $p , q \gg 0$ such
that
\begin{equation} \label{condition.eqn}
{\overline{\reg}}_H(\JJ) + \epsilon \ > \ 
\tfrac{q}{p}
\ > \ \tfrac{r_p}{p}, 
\end{equation}
where $r_p = \reg_H( \JJ^p)$. Then $\HH{1}{X}{\JJ^p
\otimes \OO_X(qH)} = 0$, and so the exact sequence
\[ 0 \lra \JJ^p \otimes \OO_X(qH) \lra \OO_X(qH) \lra
\big( \OO_X / \JJ^p \big) \otimes \OO_X(qH) \lra 0\]
together with (\ref{H0.Qp}) shows that 
\begin{equation} 
\label{ineqs} 
\begin{aligned}
\adeg_H^n \big( \JJ^p \big) \ &= \ \hh{0}{X}{\QQQ_p}
\\
	&\le \ \hh{0}{X}{(\OO_X / \JJ^p) \otimes \OO_X(qH)}
\\ &\le \ \hh{0}{X}{\OO_X(qH)}.
\end{aligned}  \end{equation}
But  Riemann-Roch implies that as a function
of $q$, 
\[ \hh{0}{X}{\OO_X(qH)} \ = \ \frac{q^n}{n!} \cdot
\deg_H(X) + O(q^{n-1}). \]
It follows from (\ref{condition.eqn}) that by taking
$p$ (and hence $q$) sufficiently large, and
$\epsilon$ sufficiently small, we can arrange that
\[ \frac{1}{p^n} \cdot \hh{0}{X}{\OO_X(qH)} \ \le \
\frac{\big( \overline{\reg}_H(\JJ) \big)^n}{n!} \cdot
\deg_H(X) + C \epsilon \]
where $C$ is a constant. The result then follows
from (\ref{ineqs}). 
\end{proof}

\begin{remark}  \textbf{(Non-Asymptotic Pathology).}
In Example \ref{Pathology.Example}, we constructed
for each $d \ge 1$  ideals $\JJ = \JJ_d$ on $\PP^3$
with fixed $s$-invariant $s_H(\JJ_d) = 2$,  but having
$d$ embedded points. This shows that one cannot bound
the arithmetic degree of an ideal in terms of its
$s$-invariant. One easily checks
that the regularity of the ideals $\JJ_d$ also goes to
infinity with $d$. So by the same token, the
regularity of a given ideal cannot be bounded in
terms of its $s$-invariant. This pathology constrasts
with results of Bayer and Mumford \cite{BM} showing
that there are (multiply exponential) bounds for the
regularity and arithmetic degree of a homogeneous
ideal in terms of the degrees of its generators. The
overall picture that seems to emerge is that the
singly exponential Bezout-type bounds appearing in
\cite{BM} are explained geometrically,
i.e. in terms of the
$s$-invariant, whereas the multiply-exponential
bounds on regularity and arithmetic degree are more
algebraic in nature. \qed
\end{remark}

\end{document}